\documentclass{article}

\usepackage{geometry}
\usepackage{lmodern} 
\usepackage{wrapfig}
\usepackage{colortbl}
\usepackage{amsmath}
\usepackage{amssymb}
\usepackage{mathrsfs}
\usepackage{makeidx}

\usepackage{multirow}
\usepackage{caption}
\usepackage{arydshln}

\newtheorem{thm}{Theorem}[section]

\usepackage{graphicx}

\usepackage{float}

\floatstyle{ruled}
\newfloat{triangle}{!ht}{tar}
\floatname{triangle}{arithmetic triangle}

\usepackage{ulem}

\usepackage{url}
\date{} 

\title{Generalized $3x+1$ Mappings : searching for cycles}

\author{Robert Tremblay\\
				Boucherville,
				Canada (Qu\'ebec),\\				
				\texttt{roberttremblay02@videotron.ca}}

\newenvironment{changemargin}[2]{%
\begin{list}{}{%
\setlength{\leftmargin}{#1}%
\setlength{\rightmargin}{#2}%
}%
\item[]}
{\end{list}}

\begin{document}


\maketitle

\begin{abstract}
We determine the conditions for the existence or not of cycles for several families of generalized $3x+1$ mappings and develop a method to find them.
\end{abstract}

\section{Introduction}

Mappings can be define on integers represented by functions such that each element of the set $\mathbb{Z}$ is connected to a single element of this set. These functions consist of two or more integer transformations on themselves. The two best known generate the original Collatz problem and the $3x+1$ problem~\cite{lagarias}. 

Let $n$ be an integer. The 3 transformations that give rise to the original Collatz problem are $2n/3$ for all integers $n=0+3q$ where $q$ is any integer, positive, negative or zero, $(4n-1)/3$ for all integers $n=1+3q$ and, $(4n+1)/3$ for all integers $n=2+3q$. The 2 transformations that apply in the $3x+1$ problem are $n/2$ and $(3n+1)/2$ respectively for even and odd integers.

The mappings generating these two problems are part of a much larger family, the generalized $3x+1$ mappings defined by Matthews~\cite{matthews_a}. The successive application of the functions that represent these mappings for any integer $n$ produces a sequence of integers called a trajectory. If we find the starting integer after $k$ operations, we have a cycle of length $k$. The cycle is then repeated to infinity. The two problems mentioned above are defined from the convergence or not of the trajectories towards the cycles. When studying the families of generalized mappings, we observe a point that appears common to all families, that the number of cycles seems limited. In the original Collatz problem the 9 known cycles are closed (there are no integers other than those included in the cycles which converges towards these cycles), which leads to the conjecture that all these other integers are in infinite trajectories (divergence). In the other known problem, only one cycle for natural numbers have found, $\langle 2,1 \rangle$, and the trajectories of all other positive integers seem to converge towards this cycle (opened cycle), leading to the famous conjecture claiming that the trajectories of all natural numbers converge towards this cycle.

Using computer programs, several cycles were determined~\cite{matthews_d} in many families of the generalized $3x+1$ mappings. Various conjectures concerning the number of these cycles, as well as the convergence or not of the trajectories, have also been stated. There are not really any methods that have been developed to determine the cycles, apart from the result of the work done by Atkin~\cite{atkin} when studying the function related to the original Collatz problem.   


In this paper we determine under which conditions a cycle can exist or not and develop a method to find them, when studying the function that generates the infinite permutations (original Collatz problem). Thereafter, we apply this method to the function related to the $3x+1$ problem and finally, to some mapping families studied by Carnielli~\cite{carnielli}. In the course of the developments, we come to a somewhat unexpected result that directly links the original Collatz problem and the $3x+1$ problem. 


\section{Infinite permutations}

Let the function $g(n)$ be defined as follows~\cite{lagarias}

\begin{equation}
g(n)=\left\lbrace  
\begin{array}{ll}
\frac{2n}{3} & \mbox{, if $n\equiv0\pmod{3}$}\\
\\
\frac{4n-1}{3} & \mbox{, if $n\equiv1\pmod{3}$}\\
\\
\frac{4n+1}{3} & \mbox{, if $n\equiv2\pmod{3}$}\\
\end{array}
\right.  
\label{CollatzOriginalTransfos}
\end{equation}

Consider the infinite permutation

\begin{equation*}
\left(
\begin{array}{cccccccccccc}
1 & 2 & 3 & 4 & 5 & 6 & 7 & 8  & 9 & \cdots & $n$    &  \cdots  \\
1 & 3 & 2 & 5 & 7 & 4 & 9 & 11 & 6 & \cdots & $g(n)$ &  \cdots
\end{array}
\right).
\end{equation*} 

The iterative application of the function to natural numbers gives rise to sequences of positive integers, called trajectories, 

\begin{equation*}
(n, g(n),g^{(2)}(n),g^{(3)}(n), \cdots, g^{(i)}(n), \cdots) ,
\end{equation*}

 with $g^{(i+1)}(n) = g\{g^{(i)}(n)\}$, $i =0, 1, 2, 3, \cdots$ and $g^{(0)}(n) = n$.

The study of the iterates of $g(n)$ is called the \textsl{original Collatz problem}. We talk about infinite permutations because when we apply the function $g$ to all positive integers a first time, we find again each of the natural numbers, but in a different order, and so on. The first transformation gives the natural integers $2 + 2q$, the second $1 + 4q$, and the third $3+4q$ where $q$ is any integer, positive or zero. 

A sequence of integers forms a loop when there exists a $k$ such that 

\begin{equation}
g^{(k)}(n) = n. 
\label{condition_cycle}
\end{equation}

If all integers in the sequence are different two by two, we have by definition a cycle of length $p = k$. Generally, we note the sequence characterizing a cycle starting with the smallest integer.

The first natural number forms a cycle noted $\langle1\rangle$. The following two numbers generate the cycle $\langle2, 3\rangle$ with a period $p = 2$. Two other cycles are known, namely

\begin{equation*}
\begin{array}{ccc}
\langle4, 5, 7, 9, 6\rangle & and & \langle44, 59, 79, 105, 70, 93, 62, 83, 111, 74, 99, 66\rangle,
\end{array}
\end{equation*}

 respectively, with the periods $p = 5$ and $p = 12$.

If we extend the problem from the set of natural numbers to the set of integers, we add the cycles~\cite{matthews_a}

\begin{equation*}
\begin{array}{cccccc}
\langle0\rangle & \langle-1\rangle & \langle-2, -3\rangle & \langle-4, -5, -7, -9, -6\rangle & and & \langle-44, -59, \cdots, -66\rangle.
\end{array}
\end{equation*}

The cycles are the same with the negative integers because the function is odd, $g(-n) = -g(n)$. In addition, the cycles are closed; there are no integers other than those included in the cycles which converges towards these cycles.

The general expression giving the result of $k$ iterations of the function $g$ on an integer $n$ is 

\begin{equation}
	g^{(k)}(n) = \lambda_{k_1,k_2}n + \rho_k(n),  
\label{g_k}
\end{equation}

where

\begin{equation}
	\lambda_{k_1,k_2} = \left(\frac{4}{3}\right)^{k_1}\left(\frac{2}{3}\right)^{k_2}  
\label{lambda_k1_k2}
\end{equation}

and 
 
\begin{equation}
	k = k_1 + k_2,
\label{k1_plus_k2}
\end{equation}

with $k_2$ the number of transformations of the form $2n/3$ and $k_1$, transformations of the other two kinds, $(4n\pm1)/3$.

Unlike parameter $\lambda_{k_1,k_2}$, $\rho_k(n)$ depend on the order of application of the transformations. Nevertheless, the maxima of this parameter are easily calculated according $k_1$.


\begin{thm}
	The absolute value of the negative or positive maximum of parameter $\rho_k(n)$ is
	
	\begin{equation}
		|\rho_{max}| = \frac{4^{k_1}-3^{k_1}}{3^{k_1}}
		\label{rho_max}
	\end{equation}
	\label{cor_rho_max}
\end{thm}

\textsl{Proof}.

We have the maxima after $k$ iterations when the $k_2$ transformations precede the $k_1$ transformations. So, for $k_1 \geq 1$,

\begin{equation*}
	|\rho_{max}| = \frac{4^{k_1-1}}{3^{k_1}} + \frac{4^{k_1-2}}{3^{k_1-1}} + \frac{4^{k_1-3}}{3^{k_1-2}} + \cdots + \frac{1}{3}
\end{equation*}

and with the same denominator,

\begin{equation*}
	|\rho_{max}| = \frac{4^{k_1-1} + 4^{k_1-2}.3 + 4^{k_1-3}.3^2 +\cdots + 3^{k_1-1}}{3^{k_1}}.
\end{equation*}

Adding and subtracting $3.4^{k_1-1}$ at the first term $4^{k_1-1}$ to the numerator, we have

\begin{equation*}
	4^{k_1-1} = 4^{k_1-1} +3.4^{k_1-1} - 3.4^{k_1-1} = 4^{k_1} - 3.4^{k_1-1}.
\end{equation*}

Adding this result at the second term to the numerator,

\begin{equation*}
\begin{array}{l}
	4^{k_1} - 3.4.4^{k_1-2} + 4^{k_1-2}.3 = 4^{k_1} - 4.(3.4^{k_1-2}) + 3.4^{k_1-2} \\
	= 4^{k_1} - 3.(3.4^{k_1-2}) = 4^{k_1} - 3^2.4^{k_1-2}.
\end{array}	
\end{equation*}

By continuing this process until the last to the numerator, leads to the expected result.  $\blacksquare$

The search for conditions that can generate a cycle leads to the analysis of the parameters $\lambda_{k_1,k_2}$ and $\rho_k(n)$ appearing in the equation \eqref{g_k}. A brief analysis of this equation when the term $\rho_k(n)$ is small in front of $\lambda_{k_1,k_2}n$, allows us to assert that $g^{(k)}(n) = n$ can be achieved for $\lambda_{k_1,k_2}$ close to 1. In the following, from results obtained by Atkin \cite{atkin}, we will show that the knowledge of the $\lambda_{k_1,k_2}$ parameter in the neighborhood of 1, determine conditions for the existence or not of a cycle. 

As that follows is very important, we recall the demonstration performed by Atkin with more details. Our final formulation will be slightly different so as to highlight the $\lambda_{k_1,k_2}$ parameter on which is based our subsequent analysis.


Consider the infinite permutation \eqref{CollatzOriginalTransfos} in the form 

\begin{equation}
	f(3n) = 2n \phantom{1234} f(3n - 2) = 4n - 3 \phantom{1234} f(3n - 1) = 4n - 1 ,
	\label{CollatzOriginalTransfos_2}
\end{equation}

applied on natural numbers, where $n$ is any integer positive. Of course, both forms lead to the same results. Mainly, the infinite permutation in the form \eqref{CollatzOriginalTransfos_2} allows us to easily built 5 other families of infinite permutations. We will use this property a little further, after the application of the theorem \ref{algorithm}.

Suppose that there is a cycle of a period $p = k$, and that $m$ is its least term. If there are $k - k_1 = k_2$ transformations of the form $f(3n) = 2n$ and $k_1$ transformations of the other two kinds with the integers $a_r$, then

\begin{equation}
	1 = \left(\frac{2}{3}\right)^{k_2}\left(\frac{4}{3}\right)^{k_1}\prod_{r = 1}^{k_1}\left(\frac{3f(a_r)}{4a_r}\right).
	\label{condition_cycle_1}
\end{equation}


With the definition \eqref{lambda_k1_k2} of $\lambda_{k_1,k_2}$ 

\begin{equation}
	1 = \lambda_{k_1,k_2}\prod_{r = 1}^{k_1}\left(\frac{3f(a_r)}{4a_r}\right).
	\label{condition_cycle_1_next}
\end{equation}

Also, for all $r$, $1 \leq r \leq k_1$, and because $f(a_r) = (4a_r \pm 1)/3$,

\begin{equation}
	-\frac{1}{4m} \leq 1 - \frac{3f(a_r)}{4a_r} \leq \frac{1}{4m},
	\label{inequality_on_f}
\end{equation}

and

\begin{equation}
	1-\frac{1}{4m} \leq \frac{3f(a_r)}{4a_r} \leq 1+\frac{1}{4m}.
	\label{inequality_on_f_next}
\end{equation}

Hence,

\begin{equation}
	(1-\frac{1}{4m})^{k_1} \leq \prod_{r = 1}^{k_1}\left(\frac{3f(a_r)}{4a_r}\right) \leq (1+\frac{1}{4m})^{k_1},
	\label{inequality_on_f_next_2}
\end{equation}

and of the equation \eqref{condition_cycle_1_next},

\begin{equation}
	(1-\frac{1}{4m})^{k_1} \leq 1/\lambda_{k_1,k_2} \leq (1+\frac{1}{4m})^{k_1}.
	\label{inequality_on_f_next_3}
\end{equation}

By applying the natural logarithm,

\begin{equation}
	{k_1}.ln(1-\frac{1}{4m}) \leq ln(1/\lambda_{k_1,k_2}) \leq {k_1}.ln(1+\frac{1}{4m}).
	\label{inequality_on_f_next_4}
\end{equation}

Now, for $0 < x < 1$ and using the Maclaurin series

\begin{equation}
	ln(1+x) = x - \frac{x^2}{2} + \frac{x^3}{3} - \frac{x^4}{4} + \cdots < x,
	\label{series_1}
\end{equation}

\begin{equation}
	\frac{1}{(1-x)} = 1 + x + x^2 + + x^3 + \cdots.
	\label{series_3}
\end{equation}

\begin{equation}
	ln(1-x) = -x - \frac{x^2}{2} - \frac{x^3}{3} - \frac{x^4}{4} - \cdots > -x - \frac{x^2}{2(1-x)},
	\label{series_2}
\end{equation}

By replacing $x$ by $1/4m$ in \eqref{series_2}, we have

\begin{equation}
	-x - \frac{x^2}{2(1-x)} = -\frac{1}{4m} - \frac{1}{4m.2(4m-1)}.
	\label{Development_4m}
\end{equation}

For $m \geq 1$,

\begin{equation}
	ln(1-\frac{1}{4m}) > -\frac{1}{4m} - \frac{1}{4m.2(4m-1)} \geq -\frac{1}{4m} - \frac{1}{4m.2.3} = - \frac{7}{24m}.
	\label{Development_4m_next}
\end{equation}

If we had chosen $m \geq 8$ (knowing that the first 7 natural numbers are already in cycles), we would have

\begin{equation}
	ln(1-\frac{1}{4m}) > -\frac{1}{4m} - \frac{1}{4m.2(4m-1)} \geq -\frac{1}{4m} - \frac{1}{4m.2.31} = - \frac{63}{248m},
	\label{Development_4m_avec_8}
\end{equation}

by putting $m = 8$ in the factor $(4m - 1)$.

Replacing $x$ by $1/4m$ in \eqref{series_1}, 

\begin{equation}
	ln(1+\frac{1}{4m}) < \frac{1}{4m} < \frac{7}{24m}.
	\label{Development_4m_next_2}
\end{equation}

Finally \eqref{inequality_on_f_next_4} becomes,

\begin{equation}
	-\frac{7}{24m} < ln(1-\frac{1}{4m}) \leq \frac{1}{k_1}ln(\frac{1}{\lambda_{k_1,k_2}}) \leq ln(1+\frac{1}{4m}) < \frac{7}{24m},
	\label{inequality_on_f_next_5}
\end{equation}

and

\begin{equation}
	\frac{1}{k_1}\left|ln(\frac{1}{\lambda_{k_1,k_2}})\right| \leq \frac{7}{24m}.
	\label{inequality_on_f_next_6}
\end{equation}
 
The condition $C$ on $m$ appears in the following inequality  

\begin{equation}
	m \leq \frac{\frac{7}{24}}{\frac{1}{k_1}\left|ln(\lambda_{k_1,k_2})\right|} = C.
	\label{condition_RT}
\end{equation}
 
For a given $k$, Atkin found

\begin{equation}
	m \leq \frac{\frac{63}{248}}{\underset{k_1}{min} \left|\frac{k}{k_1}log(1.5) - log(2)\right|} = C,
	\label{condition_Atkin}
\end{equation}

where the logarithms are in the natural base. The inequality is valid for m $\geq 8$. For a given $k (k = k_1 + k_2)$, we take the value of $k_1$ which gives the minimum of the denominator. So, we have the maximum of $C$ for this $k$.

By using $\lambda_{k_1,k_2}$ of the equation \eqref{lambda_k1_k2}, and replacing the parameter $7/24$ by $63/248$ at the numerator, the condition \eqref{condition_RT} reduces to that of Atkin \eqref{condition_Atkin}. For a given $(k_1, k_2)$, the inequalities \eqref{condition_RT} and \eqref{condition_Atkin} indicate that there can be no cycles beyond a certain $C$. The smallest integer $m$ of the cycle cannot exceed this value. These inequalities therefore impose a limit on $m$. Note that $C$ increases as $\lambda_{k_1,k_2}$ is close to 1. Conversely, $C$ decreases very rapidly as $\lambda_{k_1,k_2}$ moves away from 1. 




We will see that the analysis of $\lambda = \lambda_{k_1,k_2}$ near 1 gives not only the maxima of $C$, but also the most probable trajectories (in fact, the conditions on $k_1$ and $ k_2$ for its trajectories) to the existence of cycles.  

Here, we could directly present the theorem \ref{algorithm} and apply the resulting method which determine the values of $k_1$ and $k_2$ giving the maxima of $C$. We prefer to adopt a more inductive reasoning. We first analyze how the parameter $\lambda$ growing near 1 by adding one of the 3 transformations at a time (while remaining close to 1). Thereafter, the theorem \ref{algorithm} provides a method to find the minima of $ln(\lambda_{k_1,k_2})$ in the inequality \eqref{condition_RT} and therefore, the maxima of $C$.


First, let us a write a theorem giving the range in which the parameter $\lambda$ is located near 1. 


\begin{thm}
	The values of parameter $\lambda_k = \lambda_{k_1,k_2}$ close to 1 are between $\frac{1}{2}$ and $2$.

	\label{th_range_lambda}
\end{thm}

\textsl{Proof}.

The condition that $\lambda_{k+1} = \left(\frac{2}{3}\right)\lambda_k$ remains above 1 is that $\lambda_k > \frac{3}{2}$, otherwise $\lambda_{k+1} = \left(\frac{4}{3}\right)\lambda_k$.

If $\lambda_k = \frac{3}{2} - \epsilon$, then

\begin{equation*}
 \lambda_{k+1} = \left(\frac{4}{3}\right)\left(\frac{3}{2} - \epsilon\right) = 2 - \left(\frac{4}{3}\right)\epsilon.
\end{equation*}
	
\

The condition that $\lambda_{k+1} = \left(\frac{4}{3}\right)\lambda_k$ remains below 1 is that $\lambda_k < \frac{3}{4}$, otherwise $\lambda_{k+1} = \left(\frac{2}{3}\right)\lambda_k$.

If $\lambda_k = \frac{3}{4} + \epsilon$, then

\begin{equation*}
 \lambda_{k+1} = \left(\frac{2}{3}\right)\left(\frac{3}{4} + \epsilon\right) = \frac{1}{2} + \left(\frac{2}{3}\right)\epsilon. \phantom{12} \blacksquare
\end{equation*}  

\

Let $PP$ be $\lambda_{k_1,k_2}$ smaller than 1 ("Plus Petit que 1")  and $PG$ larger than 1 ("Plus Grand que 1"), while remaining close to 1 (by theorem~\ref{th_range_lambda}). Starting with $PP = 2/3$ and $PG = 4/3$ we have the following first results:

\begin{center}
\begin{tabular}{rl|lr}
	
	\multicolumn{2}{c|}{PP} & \multicolumn{2}{c}{PG} \\
	\hline
	\multicolumn{2}{c|}{} & \multicolumn{2}{c}{} \\
	$\frac{2}{3} = $ & \cellcolor[gray]{0.7} $\frac{2}{3}$ & \cellcolor[gray]{0.7} $\frac{4}{3}$ & $ = \frac{4}{3}$   \\ 
	\multicolumn{2}{c|}{} & \multicolumn{2}{c}{} \\
	$\frac{8}{9} = $ & \cellcolor[gray]{0.7} $\frac{2}{3}.\frac{4}{3}$ & $\frac{4}{3}.\frac{4}{3}$ & $ = \frac{16}{9}$    \\ 
	\multicolumn{2}{c|}{} & \multicolumn{2}{c}{} \\
	$\frac{16}{27} = $ & $\frac{2}{3}.\frac{4}{3}.\frac{2}{3}$ & \cellcolor[gray]{0.7} $\frac{4}{3}.\frac{4}{3}.\frac{2}{3}$ & $ = \frac{32}{27}$    \\ 
	\multicolumn{2}{c|}{} & \multicolumn{2}{c}{} \\
	$\frac{64}{81} = $ & $\frac{2}{3}.\frac{4}{3}.\frac{2}{3}.\frac{4}{3}$ & $\frac{4}{3}.\frac{4}{3}.\frac{2}{3}.\frac{4}{3}$ & $ = \frac{128}{81}$    \\ 
	\multicolumn{2}{c|}{} & \multicolumn{2}{c}{} \\
	$\frac{128}{243} = $ & $\frac{2}{3}.\frac{4}{3}.\frac{2}{3}.\frac{4}{3}.\frac{2}{3}$ & \cellcolor[gray]{0.7} $\frac{4}{3}.\frac{4}{3}.\frac{2}{3}.\frac{4}{3}.\frac{2}{3}$ & $ = \frac{256}{243}$    \\ 
	\multicolumn{2}{c|}{} & \multicolumn{2}{c}{} \\
	$\frac{512}{729} = $ & $\frac{2}{3}.\frac{4}{3}.\frac{2}{3}.\frac{4}{3}.\frac{2}{3}.\frac{4}{3}$ & $\frac{4}{3}.\frac{4}{3}.\frac{2}{3}.\frac{4}{3}.\frac{2}{3}.\frac{4}{3}$ & $ = \frac{1024}{729}$    \\ 	
	\multicolumn{2}{c|}{} & \multicolumn{2}{c}{} \\
	$\frac{2048}{2187} = $ & \cellcolor[gray]{0.7} $\frac{2}{3}.\frac{4}{3}.\frac{2}{3}.\frac{4}{3}.\frac{2}{3}.\frac{4}{3}.\frac{4}{3}$ & $\frac{4}{3}.\frac{4}{3}.\frac{2}{3}.\frac{4}{3}.\frac{2}{3}.\frac{4}{3}.\frac{4}{3}$ & $ = \frac{4096}{2187}$    \\ 		
	\multicolumn{2}{c|}{} & \multicolumn{2}{c}{} \\
	\multicolumn{2}{c|}{...} & \multicolumn{2}{c}{...} \\

\end{tabular}
\end{center}

The shaded transformations correspond to $\lambda$ giving the maxima of $C$ as we will see in he next theorem. We note these $\lambda$, $PP_{max}$ or $PG_{max}$. 

The intermediate values between two consecutive $PP_{max}$ are smaller than the two $PP_{max}$, but greater than $1/2$ (theorem \ref{th_range_lambda}). 

The intermediate values between two consecutive $PG_{max}$ are greater than the two $PG_{max}$, but smaller than $2$ (theorem \ref{th_range_lambda}).

Now, we present an method allowing us to determine the conditions on $k$, $k_1$ and $k_2$ giving the values of parameter $\lambda$ corresponding to the maxima of $C$.


\begin{thm}
	The values of $\lambda_{k_1,k_2}$ \eqref{lambda_k1_k2} calculated from the successive products of $PP$ and $PG$ correspond to the maxima of $C$ and gradually get closer to 1 with the increase of $k$. 
	\label{algorithm}
\end{thm}

\textsl{Proof}. 
	
 Let $PP = 1 - \Delta PP$ and $PG = 1 + \Delta PG$. We have the product

\begin{equation}
	PP\cdot PG=1 + \Delta PG - \Delta PP - \Delta PP\cdot \Delta PG
	\label{product}
\end{equation}	

which leads to

\begin{equation}
	1-\Delta PP < 1 + \Delta PG - \Delta PP - \Delta PP\cdot \Delta PG < 1 + \Delta PG.
	\label{inequality}
\end{equation}

$PP \cdot PG \neq 1$ and $\Delta PP \neq \Delta PG$ (except for the first two $PP$ and $PG$). Indeed, the first $\Delta PP$ and $\Delta PG$ in \textsl{base 3} are

\begin{changemargin}{2cm}{-0cm}  
$\Delta PP_1 = 3^{-1}$
\end{changemargin}

\begin{changemargin}{4cm}{-0cm}
$\Delta PG_1 = 3^{-1}$
\end{changemargin}

\begin{changemargin}{2cm}{-0cm}
$\Delta PP_2 = 3^{-2}$
\end{changemargin}

\begin{changemargin}{4cm}{-0cm}
$\Delta PG_2 = 3^{-1} - 3^{-2} - 3^{-3} = 3^{-1.535026479}$
\end{changemargin}

\begin{changemargin}{4cm}{-0cm}
$\Delta PG_3 = 3^{-3} + 3^{-4} + 3^{-5} = 3^{-2.66528248}$
\end{changemargin}

\begin{changemargin}{2cm}{-0cm}
$\Delta PP_3 = 3^{-2} - 3^{-3} - 3^{-4} + 3^{-6} + 3^{-7} = 3^{-2.508448264}$
\end{changemargin}

\begin{changemargin}{2cm}{-0cm}
$\Delta PP_4 = 3^{-4} + 3^{-6} - 3^{-7} + 3^{-8} + 3^{-9} - 3^{-11} + 3^{-12} = 3^{-3.921365509}$
\end{changemargin}

\begin{changemargin}{4cm}{-0cm}
$\Delta PG_4 = 3^{-3} + 3^{-5} - 3^{-6} - 3^{-7} - 3^{-11} + 3^{-13} - 3^{-17} = 3^{-2.945594698}$
\end{changemargin}

\begin{changemargin}{4cm}{-0cm}
$\Delta PG_5 = \cdots$
\end{changemargin}

More generally, for $t \, \epsilon \, (-1,0,+1)$,

\begin{changemargin}{2cm}{-0cm}
$\Delta PP = 3^{-r} = t_{a}3^{-a} + t_{a+1}3^{-a-1} + t_{a+2}3^{-a-2} + \cdots + t_{k_{PP}}3^{-k_{PP}}$
\end{changemargin}

\begin{changemargin}{4cm}{-0cm}
$\Delta PG = 3^{-s} = t_{b}3^{-b} + t_{b+1}3^{-b-1} + t_{b+2}3^{-b-2} + \cdots + t_{k_{PG}}3^{-k_{PG}}$
\end{changemargin}

with $t_a=t_b=1$ , $t_{k_{PP}} \neq 0$ and $t_{k_{PG}} \neq 0$.

The exponent, in absolute value, of the last term (the smallest) of each $\Delta PP$ (or $\Delta PG$) is equal to $k_{PP}$ (or $k_{PG}$), so the number of transformations $k=k_1+k_2$.

Thus, the successive products of $PP \cdot PG$ give values which approach more and more of $1$ without ever reaching it and according to equation \eqref{condition_RT}, we find the maxima of $C$, which leads to the following algorithm.




$\blacksquare$

\subsection*{Algorithm}

Start with $PP = 2/3$ and $(k_1, k_2) = (0, 1)$; $PG = 4/3$ and $(k_1, k_2) = (1, 0)$.

The first product is $PP \cdot PG = 2/3 \cdot 4/3 = 8/9$, and $(k_1, k_2) = (0 + 1, 1 + 0) = (1, 1)$. This operation determines a new $PP$, $PP = 8/9$. 

The product of this new $PP$ with $PG$ gives a new $PG$, $PP \cdot PG = 8/9 \cdot 4/3 = 32/27$, and $(k_1, k_2) = (1 + 1, 1 + 0) = (2, 1)$. The new $PG$ is $PG = 32/27$. Then, $PP = 8/9$ and $PG = 32/27$. The product $PP \cdot PG = 256/243$ identify a new $PG$, and $(k_1, k_2) = (1 + 2, 1 + 1) = (3, 2)$. 

By repeating this process we get the pairs $(k_1, k_2)$ that give the $PP$ and the $PG$ corresponding to the maxima of $C$.  

Define a node as the set of consecutive maxima ($PP_{max}$ or $PG_{max}$). Let $i$ the parameter identifying the sequence of these sets (primary or main nodes) and $j$ each element of its sets (secondary nodes). Then the notation $N_{i,j}$ represent all the nodes with $N_{1,1}$ identifying both, the first $PP_{max} = 2/3$ and the first $PG_{max} = 4/3$.

\subsection*{Results}

The results for the first nine nodes are presented in the table~\ref{NodesInfinitePermutation}.  

The pair of integers $(k_1,k_2)$ obtained by the preceding algorithm determines the maxima of $C$. Let $m$ be the least integer of a trajectory generated by the transformation~\eqref{CollatzOriginalTransfos} or \eqref{CollatzOriginalTransfos_2} with a given $(k_1,k_2)$. Then, the only possible cycles are those for $m \leq C$. As we will see, these combinations $(k_1,k_2)$ also seem to determine the conditions for the most probable trajectories for the existence of cycles.


The $(k_1,k_2)$ of the 4 known cycles for the natural numbers are exactly equal to those of the nodes $N_{1,1}, N_{2,1}, N_{3,2}$ and $N_{4,2}$ of the table for the lenghts $p = k = 1,2,5$ and $12$.


In the table~\ref{ExampleTrajectories_K53} we give some examples of trajectories for the node $N_{6,1}$ with $k = 53$ and $(k_1,k_2) = (31, 22)$ ($\lambda$ is close to 1). We present two trajectories for this node and each couple $(30, 23)$, $(32, 21)$ and $(33, 20)$ around of the node $N_{6,1}$. 

In the next table~\ref{ExampleTrajectoriesClosedOne} we have chosen four other cases of $\lambda$ close to 1 derived of the algorithm, so $k = 17$, $29$, $41$ and $94$ corresponding respectively to the nodes $N_{5,1}$, $N_{5,2}$, $N_{5,3}$ and $N_{7,1}$. 

\textsl{Then, the cycles seem more probable for $\lambda$ close to 1 and this probability decreases very rapidly as $\lambda$ moves away from 1.}



There are some interesting families of infinite permutations that are built from the permutation of function~\eqref{CollatzOriginalTransfos_2}. Then, we have six permutations in starting with this function, 

\begin{equation}
	f(3n) = 2n \phantom{1234} f(3n - 2) = 4n - 3 \phantom{1234} f(3n - 1) = 4n - 1,
	\label{CollatzOriginalTransfos_3}
\end{equation}

\begin{equation}
	f(3n) = 2n \phantom{1234} f(3n - 2) = 4n - 1 \phantom{1234} f(3n - 1) = 4n - 3,
	\label{CollatzOriginalTransfos_4}
\end{equation}

\begin{equation}
	f(3n) = 4n - 3 \phantom{1234} f(3n - 2) = 4n - 1 \phantom{1234} f(3n - 1) = 2n,
	\label{CollatzOriginalTransfos_5}
\end{equation}

\begin{equation}
	f(3n) = 4n - 3 \phantom{1234} f(3n - 2) = 2n \phantom{1234} f(3n - 1) = 4n - 1,
	\label{CollatzOriginalTransfos_6}
\end{equation}

\begin{equation*}
	\cdots
\end{equation*}

The third function generates, among others, a cycle of period $k = 94$ for the smallest term $m = 144$. This period is associated with the first secondary node of the node 7, namely $N_{7,1}$. 

In the other form the function is

\begin{equation}
h(n)=\left\lbrace  
\begin{array}{ll}
\frac{4n-9}{3} & \mbox{, if $n\equiv0\pmod{3}$}\\
\\
\frac{4n+5}{3} & \mbox{, if $n\equiv1\pmod{3}$}\\
\\
\frac{2n+2}{3} & \mbox{, if $n\equiv2\pmod{3}$}\\
\end{array}
\right.  
\label{CollatzOriginalTransfos_7}
\end{equation}

or more simply, after a first application of $h(n)$ on the natural numbers

\begin{equation*}
\left(
\begin{array}{cccccccccccc}
1 & 2 & 3 & 4 & 5 & 6 & 7 & 8  & 9 & \cdots & $n$    &  \cdots  \\
3 & 2 & 1 & 7 & 4 & 5 & 11 & 6 & 9 & \cdots & $h(n)$ &  \cdots
\end{array}
\right).
\end{equation*} 

Also, there is a cycle for $k = 6$ ($k_1 = k_2 = 3$), namely $\langle 4, 7, 11, 8, 6, 5 \rangle$. This case does not appear in table~\ref{NodesInfinitePermutation}, but $\lambda_{k_1, k_2} = 0.70233196159122$, which is close to 1. 



\section{Problem 3x + 1}

Let us apply the search method of the cycles as developed in the iterative application of the function $g(n)$ to another similar function.

Let the function $T(n)$ be defined as follow~\cite{lagarias}

\begin{equation}
T(n)=\left\lbrace  
\begin{array}{ll}
\frac{n}{2} & \mbox{, si $n\equiv0\pmod{2}$}\\
\\
\frac{3n+1}{2} & \mbox{, si $n\equiv1\pmod{2}$}.\\
\end{array}
\right.  
\label{Trois_x_plus_1_Transfos}
\end{equation}


The iterative application of $T(n)$ on the integers generate the different trajectories. The result is the same as using the function

\begin{equation}
	f(2n) = n \phantom{1234} f(2n + 1) = 3n + 2.
	\label{transfos_4}
\end{equation}


Suppose that there is a cycle of a period $p = k$, and that $m$ is its least term. If there are $k - k_1 = k_2$ transformations of the form $f(2n) = n$ and $k_1$ transformations of the other kind with the integers $a_r$, then

\begin{equation}
	1 = \left(\frac{1}{2}\right)^{k_2}\left(\frac{3}{2}\right)^{k_1}\prod_{r = 1}^{k_1}\left(\frac{2f(a_r)}{3a_r}\right).
	\label{condition_cycle_2}
\end{equation}


For the positive integers, 

\begin{equation}
	1 - \frac{2f(a_r)}{3a_r} \geq -\frac{1}{3m}.
	\label{inequality_on_f_2}
\end{equation}

For the negative integers,

\begin{equation}
	1 - \frac{2f(a_r)}{3a_r} \leq \frac{1}{3|m|}.
	\label{inequality_on_f_3}
\end{equation}

The condition to the existence of a cycle is given by the expression

\begin{equation}
	|m| \leq \frac{\frac{5}{12}}{\frac{1}{k_1}\left|ln(\lambda_{k_1,k_2})\right|} = C.
	\label{condition_RT_2}
\end{equation}


This inequality is valid for $|m| \geq 1$. $\lambda_{k_1,k_2}$ is given by

\begin{equation}
	\lambda_{k_1,k_2} = \left(\frac{3}{2}\right)^{k_1}\left(\frac{1}{2}\right)^{k_2}.  
\end{equation}

The values of parameter $\lambda_k$ close to 1 are between 1/3 and 3.
	
The algorithm developed in the previous section generates the table~\ref{NodesProblem_3xPlus1}. The pair of integers $(k_1,k_2)$ determines the maxima of $C$. Let $m$ be the least integer of a trajectory generated by the transformations~\eqref{Trois_x_plus_1_Transfos} or \eqref{transfos_4} with a given $(k_1,k_2)$. Then, the only possible cycles are those for $|m| \leq C$. These combinations $(k_1,k_2)$ also seem to determine the conditions for the most probable trajectories for the existence of cycles.

It is interesting to note that all $PP$ and $PG$ obtained by the successive products of $PP \cdot PG$ (except $PP = 1/2$) are the reciprocals of those obtained in the infinite permutations. $PP = 1/2$, $PG = 3/2$, $PP=3/4$, $PP=9/8$, $\cdots$, in the $3x + 1$ problem and $PP = 2/3$, $PG = 4/3$, $PP=8/9$, $\cdots$ in the problem of infinite permutations. The distribution of the primary and secondary nodes is then identical. For example, in the table~\ref{NodesInfinitePermutation} (infinite permutations), node $9$ contains $23$ secondary nodes, exactly like node $10$ in the table~\ref{NodesProblem_3xPlus1}. Therefore, there is a direct link between two problems that appear when looking for the maxima of $C$, or if we want, when searching for the most probable trajectories to existence for a cycle.    
	
The general expression giving the result of $k$ iterations of the function $T$ on an integer $n$ is 

\begin{equation}
	T^{(k)}(n) = \lambda_{k_1,k_2}n + \delta_k(n).  
\label{T_k}
\end{equation}

Because $\delta_k(n)$ is always positive, the possible cycles for the positive integers are values $\lambda = PP < 1$, and for the negative integers we have the possible cycles for $\lambda = PG > 1$, with $k_1$ and $k_2$ giving $\lambda$ close to 1.

For the positive integers we have the cycle $\langle1, 2\rangle$ with the length $k = 2$ and $PP = 0.75$ corresponding to the node $N_{2,1}$ in the table.

For the zero and negative integers we have the cycles $\langle0\rangle$, $\langle-1\rangle$, $\langle-5, -7, -10\rangle$ and the long cycle~\cite{matthews_a} 

\begin{equation*}
	\langle-17, -25, -37, -55, -82, -41, -61, -91, -136, -68, -34\rangle
\end{equation*}

with lengths $k = 1$, $k = 3$ and $k = 11$.

These last values of $k$ are in table for nodes $N_{1,1}$ (PP = 0.5 and PG = 1.5), $N_{3,1}$ and $N_{5,1}$.


In the table~\ref{ExampleTrajectories_positiveIntegers} and \ref{ExampleTrajectories_negativeIntegers} we give some examples of trajectories respectively for positive integers and negative integers. \textsl{Also, the cycles seem more probable for $\lambda$ close to $1$ and this probability decreases very rapidly as $\lambda$ moves away from $1$.}

\section{Generalized 3x + 1 mappings}

Defining the \textsl{generalized Collatz mapping} or \textsl{generalized $3x + 1$ mapping}~\cite{matthews_a}

\begin{equation}	
	T(x) = \frac{m_ix - r_i}{d} \phantom{123} \mbox{, if  $x\equiv i\pmod{d},$}  
\end{equation}

 $d \geq 2$ be a positive integer and $m_0,\cdots,m_{d-1}$ be non-zero integers. Also for $i = 0,\cdots,d - 1$, let \mbox{$r_i \phantom{1} \epsilon \phantom{1} \mathbb{Z}$} satisfy $r_i\equiv im_i\pmod{d}$.

The original Collatz mapping corresponds to parameter choices $d = 3$, $m_0 = 2$, $m_1 = m_2 = 4$, $r_0 = 0$, $r_1 = 1$ and $r_2 = -1$. The $3x + 1$ mapping corresponds to the choices $d = 2$, $m_0 = 1$, $m_1 = 3$, $r_0 = 0$ and $r_1 = -1$.

Carnielli~\cite{carnielli,matthews_b} has proposed two natural generalizations of Collatz Problem which are the special cases of a generalized $3x + 1$ mapping. Let $m_0 = 1$, $r_0 = 0$ for $i = 0$, and $m_i = d + 1$ and $r_i = -(d - i)$ for $i \geq 1$

\begin{equation}
T_d(x)=\left\lbrace  
\begin{array}{ll}
\frac{x}{d} & \mbox{, if $x\equiv0\pmod{d}$}\\
\\
\frac{((d+1)x+d-i)}{d} & \mbox{, if $x\equiv i\pmod{d}, 1 \leq i \leq d-1$.}\\
\\
\end{array}
\right.  
\label{CarnielliTransfos_1}
\end{equation}

Let $m_0 = 1$, $r_0 = 0$ for $i = 0$, and $m_i = d + 1$ and $r_i = i$ for $i \neq 1$

\begin{equation}
L_d(x)=\left\lbrace  
\begin{array}{ll}
\frac{x}{d} & \mbox{, if $x\equiv0\pmod{d}$}\\
\\
\frac{((d+1)x-i)}{d} & \mbox{, if $x\equiv i\pmod{d}, -d/2 < i \leq d/2, i \neq 0$.}\\
\\
\end{array}
\right.  
\label{CarnielliTransfos_2}
\end{equation}

Keith Matthews developed the last transformation which is a generalization of the mapping of Lu Pei~\cite{matthews_c} with $d = 3$. 

The conditions for the existence of a cycle in the generalization~\eqref{CarnielliTransfos_1} or \eqref{CarnielliTransfos_2} are similar to the conditions~\eqref{condition_cycle_1} of the original Collatz problem or the condition~\eqref{condition_cycle_2} of the $3x + 1$ problem, then

\begin{equation}
	1 = \left(\frac{1}{d}\right)^{k_1}\left(\frac{d+1}{d}\right)^{k_2}\left(\frac{d+1}{d}\right)^{k_3}\left(\cdots\right)\prod_{r=1}^{k_2}\left(\frac{d}{(d+1)}.\frac{f(a_r)}{a_r}\right)\prod_{s=1}^{k_3}\left(\frac{d}{(d+1)}.\frac{f(a_s)}{a_r}\right)\left(\cdots\right).
	\label{condition_cycle_3}
\end{equation}



$k_1, k_2, k_3, \dots$ correspond to the number of transformations for $i = 0, 1, 2, \cdots$.

Also, we develop the condition on the least term $m$ of cycle and find 

\begin{equation}
	|m| \leq \frac{par}{\frac{1}{k_2+k_3+\cdots}\left|log(\lambda_{k_1,k_2, \cdots})\right|} = C,
	\label{condition_RT_3}
\end{equation}

where $par$ is a parameter and $\lambda_{k_1,k_2, \cdots}$ is given by

\begin{equation}
	\lambda_{k_1,k_2,k_3,\cdots} = \left(\frac{1}{d}\right)^{k_1}\left(\frac{d+1}{d}\right)^{k_2+k_3+\cdots}.  
\label{lambda_Carnielli}
\end{equation}

We have the maxima for $C$ when $\lambda$ is close to 1 and we can apply the algorithm developed in this paper.

Carnielli has produced two tables for $2 \leq d \leq 150$ and $|x| \leq 250,000$ giving the least term of cycles and the cycle lengths. 

We can verify a very important result : \textsl{all cycle lengths correspond to some values of $\lambda$ close to 1 and we find them with our algorithm.}  


We can apply our algorithm because there are two different terms $m_i/d$ in equation~\eqref{lambda_Carnielli}. It is more complicated when we have three or more different terms $m_i/d$ giving $\lambda$.

For example, let the generalized mapping~\cite{matthews_a} be 

\begin{equation}
T(x)=\left\lbrace  
\begin{array}{ll}
\frac{x}{4} & \mbox{, if $x\equiv0\pmod{4}$}\\
\\
\frac{3x-3}{4} & \mbox{, if $n\equiv1\pmod{4}$}\\
\\
\frac{5x-2}{4} & \mbox{, if $n\equiv2\pmod{4}$}\\
\\
\frac{17x-3}{4} & \mbox{, if $n\equiv3\pmod{4}$}.\\
\end{array}
\right.  
\label{GeneralizedTransfos_ex1}
\end{equation}

Matthews has found $17$ cycles (lengths in parentheses), starting at values 0(1), -3(1), 2(1), 3(2), 6(1747), -18(2), -46(34), -122(8), -330(4), -117(4), -137(4), -186(4), -513(1426), -261(4), -333(4), 5127(14), -5205(60). 

It is probably possible to prove that the condition $C$ on the least term of a cycle is proportional to the inverse of $ln(\lambda)$ with

\begin{equation}
	\lambda_{k_1,k_2,k_3,k_4} = \left(\frac{m_0}{d}\right)^{k_1}\left(\frac{m_1}{d}\right)^{k_2}\left(\frac{m_2}{d}\right)^{k_3}\left(\frac{m_3}{d}\right)^{k_4}. 
\end{equation}

The maxima of $C$ stands for $\lambda$ close to 1.

The six cycles of period $p=4$ correspond to $k_1=k_2=k_3=k_4=1$ and 

\begin{equation*}
	\lambda = (1/4)(3/4)(5/4)(17/4) = 0.99609375. 
\end{equation*}

For example, the cycle starting with $-330$ is $\langle -330,-413,-1756,-439 \rangle$, where the order of the transformations is $k_3, k_4, k_1, k_2$.



The cycle of period $p=8$ has $k_1=k_2=k_3=k_4=2$ and $\lambda=0.992202759$. The cycle of period $p=14$ has $k_1=2$, $k_2=7$, $k_3=2$, $k_4=3$ and $\lambda=1.000681054$. The cycle of period $p=1747$ corresponds to $k_1=432$, $k_2=434$, $k_3=450$, $k_4=431$, and $\lambda=1.354586564$.

\



\section{Conclusion}

For several families of the generalized $3x + 1$ mappings which include 2 different terms $m_i/d$, we have been able to construct a quantity $C$ such that the least integer $m$ of a trajectory generating a cycle is subjected to the condition $|m| \leq C$. These families have 2 types of transformations characterized by the quantities $(k_1, k_2)$. $k_1$ is the number of transformations of one kind and $k_2$ the set of transformations of the other kinds, with $k = k_1 + k_2$ the total number of transformations. We then built an algorithm giving the values $(k_1, k_2)$ corresponding to the maxima of $C$. Moreover the conditions $(k_1, k_2)$ seem to be those determining the possible cycles.

So, we have developed a method to determine under which conditions there is or not cycle for several families of the generalized $3x + 1$ mappings.

At the end of the paper, we discussed a case with 4 terms $m_i/d$. As we have noted, the 17 known cycles are found for the combinations $(k_1, k_2, k_3, k_4)$ generating a $\lambda$ close to 1. However, we have not demonstrated that there a quantity $C$ as in the previous cases and, if we can built an algorithm determining the maxima of $C$.

An interesting question which remains outstanding and that we raised in the introduction, is the number of cycles is limited or not ? 


\bibliographystyle{unsrt}



\begin{table}   
\begin{center}
\begin{tabular}{|c|r|l|l|r|r|r|r|}
	\hline
	\multicolumn{4}{|l|}{Main nodes} & \multicolumn{4}{r|}{} \\
	\cline{1-4}
	\multirow{2}*{} & \multicolumn{3}{|l|}{Secondary nodes} & \multicolumn{4}{r|}{} \\
	\cline{2-8} 
	\ & & \multicolumn{1}{|l|}{PP} & \multicolumn{1}{|l|}{PG} & $k_1$ & $k_2$ & k & ln (C) \\
	\hline
	\multicolumn{8}{|c|}{} \\
	\hline
	
	1 & 1 & \multicolumn{2}{|l|}{0.66666666666667} & 0 & 1 & 1 &  \\
	\hline
	1 & 1 & & 1.33333333333333 & 1 & 0 & 1 &  \\ 
	\hline
	\multicolumn{8}{|c|}{} \\
	\hline
	
	2 & 1 & \multicolumn{2}{|l|}{0.88888888888889} & 1 & 1 & 2 & 0.9067673 \\ 
	\hline		
	\multicolumn{8}{|c|}{} \\
	\hline	
	
	\multirow{2}*{3} & 1 & & 1.18518518518519 & 2 & 1 & 3 & 1.2335544  \\ 
	\cline{2-8}		
	 & 2 & & 1.05349794238683 & 3 & 2 & 5 & 2.8207519  \\ 
	\hline			
	\multicolumn{8}{|c|}{} \\
	\hline
	
	\multirow{2}*{4} & 1 & \multicolumn{2}{|l|}{0.93644261545496} & 4 & 3 & 7 & 2.8773089  \\ 
	\cline{2-8}		
	 & 2 & \multicolumn{2}{|l|}{0.98654036854514} & 7 & 5 & 12 & 5.0150589  \\ 
	\hline			
	\multicolumn{8}{|c|}{} \\
	\hline			
	
	\multirow{3}*{5} & 1 & & 1.03931824834386 & 10 & 7 & 17 & 4.3258524  \\ 
	\cline{2-8}		
	 & 2 & & 1.02532940775684 & 17 & 12 & 29 & 5.2893919  \\ 
	\cline{2-8}			
	 & 3 & & 1.01152885180861 & 24 & 17 & 41 & 6.4145496  \\ 
	\hline		
	\multicolumn{8}{|c|}{} \\
	\hline

	6 & 1 & \multicolumn{2}{|l|}{0.99791404625731} & 31 & 22 & 53 & 8.3733287  \\ 
	\hline		
	\multicolumn{8}{|c|}{} \\
	\hline

	\multirow{5}*{7} & 1 & & 1.00941884941434 & 55 & 39 & 94 & 7.4449229  \\ 
	\cline{2-8}		
	 & 2 & & 1.00731324838746 & 86 & 61 & 147 & 8.1439169  \\ 
	\cline{2-8}			
	 & 3 & & 1.00521203954693 & 117 & 83 & 200 & 8.7894147  \\ 
	\cline{2-8}			
	 & 4 & & 1.00311521373084 & 148 & 105 & 253 & 9.5380817  \\ 
	\cline{2-8}			
	 & 5 & & 1.00102276179641 & 179 & 127 & 306 & 10.841002  \\ 
	\hline		
	\multicolumn{8}{|c|}{} \\
	\hline
	
	\multirow{2}*{6} & 1 & \multicolumn{2}{|l|}{0.99893467461992} & 210 & 149 & 359 & 10.958906  \\ 
	\cline{2-8}		
	 & 2 & \multicolumn{2}{|l|}{0.99995634684222} & 389 & 276 & 665 & 14.7706488  \\ 
	\hline			
	\multicolumn{8}{|c|}{} \\
	\hline

	\multirow{5}*{9} & 1 & & 1.00097906399185 & 568 & 403 & 971 & 12.0393806  \\ 
	\cline{2-8}		
	 & 2 & & 1.00093536809484 & 957 & 679 & 1,636 & 12.6066976  \\ 
	\cline{2-8}			
	 & \dots & &  &  &  &  &   \\ 
	\cline{2-8}			
	 & 22 & & 1.00006185061131 & 8,737 & 6,199 & 14,936 & 17.533998  \\ 
	\cline{2-8}			
	 & 23 & & 1.00001819475356 & 9,126 & 6,475 & 15,601 & 18.801125  \\ 
	\hline		
				
\end{tabular}
\end{center}
\caption{Nodes - Infinite Permutations}
\label{NodesInfinitePermutation}
\end{table}


\begin{table}
\begin{center}
\begin{tabular}{|l|c|l|l|}
	\hline
	$k$ & $(k_1,k_2)$ & $\lambda_{k_1,k_2}$ & trajectories \\
	\hline
	53 & $(30,23)$ & 0.49895703128654 & $(160, 213, \cdots, 77)$ \\
  & & & $(312, 208, \cdots, 161)$ \\
	\hline
	 & \cellcolor[gray]{0.9}$(31,22)$ & 0.997914046257308 & $(225, 150, \cdots, 224)$ \\
  & & & $(326, 435, \cdots, 325)$ \\
	\hline
	 & $(32,21)$ & 1.99582809251462 & $(84, 56, \cdots, 163)$ \\
  & & & $(320, 427, \cdots, 642)$ \\
	\hline	
	 & $(33,20)$ & 3.99165618502923 & $(56, 75, \cdots, 217)$ \\
  & & & $(243, 162, \cdots, 1,007)$ \\
	\hline

\end{tabular}
\end{center}
\caption{Examples - Trajectories closed $(k_1, k_2) = (31, 22)$ - Infinite permutations}
\label{ExampleTrajectories_K53}
\end{table}


\begin{table}
\begin{center}
\begin{tabular}{|l|c|l|l|}
	\hline
	$k$ & $(k_1,k_2)$ & $\lambda_{k_1,k_2}$ & trajectories \\
	\hline
	17 & $(10,7)$ & 1.03931824834385 & $(36, 24, \cdots, 37)$ \\
  & & & $(46, 61, \cdots, 47)$ \\
	\hline
	29 & $(17,12)$ & 1.02532940775684 & $(78, 52, \cdots, 77)$ \\
  & & & $(88, 117, \cdots, 87)$ \\
	\hline
	41 & $(24,17)$ & 1.01152885180861 & $(50, 67, \cdots, 49)$ \\
  & & & $(448, 597, \cdots, 449)$ \\
	\hline	
	94 & $(55,39)$ & 1.00941884941434 & $(332, 443, \cdots, 336)$ \\
  & & & $(742, 989, \cdots, 745)$ \\
	\hline

\end{tabular}
\end{center}
\caption{Examples - Trajectories closed to 1 - Infinite permutations}
\label{ExampleTrajectoriesClosedOne}
\end{table}


\begin{table}  
\begin{center}
\begin{tabular}{|c|r|l|l|r|r|r|r|}
	\hline
	\multicolumn{4}{|l|}{Main nodes} & \multicolumn{4}{r|}{} \\
	\cline{1-4}
	\multirow{2}*{} & \multicolumn{3}{|l|}{Secondary nodes} & \multicolumn{4}{r|}{} \\
	\cline{2-8} 
	\ & & \multicolumn{1}{|l|}{PP} & \multicolumn{1}{|l|}{PG} & $k_1$ & $k_2$ & k & ln (C) \\
	\hline
	\multicolumn{8}{|c|}{} \\
	\hline
	
	1 & 1 & \multicolumn{2}{|l|}{0.50000000000000} & 0 & 1 & 1 &   \\ 
	\hline
	1 & 1 & & 1.500000000000000 & 1 & 0 & 1 &   \\ 
	\hline
	\multicolumn{8}{|c|}{} \\
	\hline
	
	2 & 1 & \multicolumn{2}{|l|}{0.75000000000000} & 1 & 1 & 2 & 0.3704306  \\ 
	\hline		
	\multicolumn{6}{|c|}{} \\
	\hline	
	
	3 & 1 & & 1.12500000000000 & 2 & 1 & 3 & 1.9565895  \\ 
	\hline		
	\multicolumn{8}{|c|}{} \\
	\hline
	
	\multirow{2}*{4} & 1 & \multicolumn{2}{|l|}{0.84375000000000} & 3 & 2 & 5 & 1.9956945  \\ 
	\cline{2-8}		
	 & 2 & \multicolumn{2}{|l|}{0.94921875000000} & 5 & 3 & 8 & 3.6882524  \\ 
	\hline			
	\multicolumn{8}{|c|}{} \\
	\hline			

	\multirow{2}*{5} & 1 & & 1.06787109375000 & 7 & 4 & 11 & 3.7935996  \\ 
	\cline{2-8}		
	 & 2 & & 1.01364326477050 & 12 & 7 & 19 & 5.9107304  \\ 
	\hline			
	\multicolumn{8}{|c|}{} \\
	\hline	
	
	\multirow{3}*{6} & 1 & \multicolumn{2}{|l|}{0.96216919273138} & 17 & 10 & 27 & 5.2131556  \\ 
	\cline{2-8}		
	 & 2 &\multicolumn{2}{|l|}{0.97529632178184} & 29 & 17 & 46 & 6.1801493  \\ 
	\cline{2-8}			
	 & 3 & \multicolumn{2}{|l|}{0.98860254772961} & 41 & 24 & 65 & 7.3067428  \\ 
	\hline		
	\multicolumn{8}{|c|}{} \\
	\hline

	7 & 1 & & 1.00209031404109 & 53 & 31 & 84 & 9.2663084  \\ 
	\hline		
	\multicolumn{8}{|c|}{} \\
	\hline

	\multirow{5}*{8} & 1 & \multicolumn{2}{|l|}{0.99066903751619} & 94 & 55 & 149 & 8.3375594  \\ 
	\cline{2-8}		
	 & 2 & \multicolumn{2}{|l|}{0.99273984691538} & 147 & 86 & 233 & 9.0366771  \\ 
	\cline{2-8}			
	 & 3 & \multicolumn{2}{|l|}{0.99481498495653} & 200 & 117 & 317 & 9.6822330  \\ 
	\cline{2-8}			
	 & 4 & \multicolumn{2}{|l|}{0.99689446068787} & 253 & 148 & 401 & 10.4309339  \\ 
	\cline{2-8}			
	 & 5 & \multicolumn{2}{|l|}{0.99897828317652} & 306 & 179 & 485 & 11.7338762  \\ 
	\hline		
	\multicolumn{8}{|c|}{} \\
	\hline
	
	\multirow{2}*{9} & 1 & & 1.00106646150859 & 359 & 210 & 569 & 11.8517958  \\ 
	\cline{2-8}		
	 & 2 & & 1.00004365506344 & 665 & 389 & 1,054 & 15.6635314  \\ 
	\hline			
	\multicolumn{8}{|c|}{} \\
	\hline

	\multirow{5}*{8} & 1 & \multicolumn{2}{|l|}{0.99902189363685} & 971 & 568 & 1,539 & 12.9322606  \\ 
	\cline{2-8}		
	 & 2 & \multicolumn{2}{|l|}{0.99906550600100} & 1,636 & 957 & 2,593 & 13.4995787  \\ 
	\cline{2-8}			
	 & \dots & &  &  &  &  &   \\ 
	\cline{2-8}			
	 & 22 & \multicolumn{2}{|l|}{0.99993815321363} & 14,936 & 8,737 & 23,673 & 18.426880  \\ 
	\cline{2-8}			
	 & 23 & \multicolumn{2}{|l|}{0.99998180557715} & 15,601 & 9,126 & 24,727 & 19.694008  \\ 
	\hline
				
\end{tabular}
\end{center}
\caption{Nodes - 3x + 1}
\label{NodesProblem_3xPlus1}
\end{table}


\begin{table}
\begin{center}
\begin{tabular}{|l|c|l|l|}
	\hline
	$k$ & $(k_1,k_2)$ & $\lambda_{k_1,k_2}$ & trajectories \\
	\hline
	27 & $(16,11)$ & 0.320723064243793 & $(159, 239, \cdots, 53)$ \\
  & & & $(239, 359, \cdots, 80)$ \\
	\hline
	 & \cellcolor[gray]{0.9}$(17,10)$ & 0.96216919273138 & $(166, 83, \cdots, 167)$ \\
  & & & $(250, 125, \cdots, 251)$ \\
	\hline
	 & $(18,9)$ & 2.88650757819414 & $(54, 27, \cdots, 167)$ \\
  & & & $(82, 41, \cdots, 251)$ \\
	\hline	
	 & $(19,8)$ & 8.65952273458242 & $(27, 41, \cdots, 251)$ \\
  & & & $(31, 47, \cdots, 283)$ \\
	\hline	
	46 & \cellcolor[gray]{0.9}$(29,17)$ & 0.97529632178194 & $(91, 137, \cdots, 92)$ \\
  & & & $(121, 182, \cdots, 122)$ \\
	\hline
	65 & \cellcolor[gray]{0.9}$(41,24)$ & 0.988602547729613 & $(73, 110, \cdots, 80)$ \\
  & & & $(231, 347, \cdots, 244)$ \\
	\hline

\end{tabular}
\end{center}
\caption{Examples - Trajectories positive integers - 3x + 1}
\label{ExampleTrajectories_positiveIntegers}
\end{table}


\begin{table}
\begin{center}
\begin{tabular}{|l|c|l|l|}
	\hline
	$k$ & $(k_1,k_2)$ & $\lambda_{k_1,k_2}$ & trajectories \\
	\hline
	11 & $(5,6)$ & 0.11865234375 & $(-63, -94, \cdots, -7)$ \\
  & & & $(-69, -103, \cdots, -8)$ \\
	\hline
	 & $(6,5)$ & 0.35595703125 & $(-49, -73, \cdots, -17)$ \\
  & & & $(-81, -121, \cdots, -28)$ \\
	\hline
	 & \cellcolor[gray]{0.9}$(7,4)$ & 1.06787109375 & $(-42, -21, \cdots, -41)$ \\
  & & & $(-57, -85, \cdots, -59)$ \\
	\hline	
	 & $(8,3)$ & 3.20361328125 & $(-21, -31, \cdots, -61)$ \\
  & & & $(-145, -217, \cdots, -460)$ \\
	\hline	
	19 & \cellcolor[gray]{0.9}$(12,7)$ & 1.0136932647051 & $(-65, -97, \cdots, -64)$ \\
  & & & $(-87, -130, \cdots, -86)$ \\
	\hline		
	
\end{tabular}
\end{center}
\caption{Examples - Trajectories negative integers - 3x + 1}
\label{ExampleTrajectories_negativeIntegers}
\end{table}

\end{document}